This article has been accepted for publication in IEEE Transactions on Automatic Control. This is the author's version which has not been fully edited and content may change prior to final publication. Citation information: DOI 10.1109/TAC.2022.3214795

# Fixed-Time Stable Proximal Dynamical System for Solving MVIPs

Kunal Garg, Mayank Baranwal, Rohit Gupta and Mouhacine Benosman

*Abstract*—**In this paper, a novel modified proximal dynamical system is proposed to compute the solution of a mixed variational inequality problem (MVIP) within a fixed time, where the time of convergence is finite and is uniformly bounded for all initial conditions. Under the assumptions of strong monotonicity and Lipschitz continuity, it is shown that a solution of the modified proximal dynamical system exists, is uniquely determined, and converges to the unique solution of the associated MVIP within a fixed time. Furthermore, the fixed-time stability of the modified projected dynamical system continues to hold, even if the assumption of strong monotonicity is relaxed to that of strong pseudomonotonicity. Finally, it is shown that the solution obtained using the forward-Euler discretization of the proposed modified proximal dynamical system converges to an arbitrarily small neighborhood of the solution of the associated MVIP within a fixed number of time steps, independent of the initial conditions.**

*Index Terms*—**Mixed variational inequality problem; Proximal dynamical system; Fixed-Time stability; Discretization**

## I. Introduction

Mixed variational inequality problems (MVIPs) have numerous applications in optimization (see, e.g., [1], [2]), game theory (see, e.g., [3], [4]), control theory (see, e.g., [5], [6]), and other related areas (see, e.g., [7]). In the literature, both discrete-time gradient-based methods, and continuous-time gradient flow-based approaches have been proposed for solving MVIPs. The focus of this paper is on designing continuous-time dynamical systems such that their solutions converge to the solution of the MVIP, described in (1) (see Section II), in a fixed time, starting from any given initial condition.

In recent years, the use of dynamical systems has emerged as a viable alternative for solving MVIPs with a particular focus on optimization problems (see, e.g., [3], [8], [9], [10]). This viewpoint allows tools from Lyapunov theory to be employed for the stability analysis of the equilibrium points of the underlying dynamical systems. Under the assumptions of monotonicity and strong monotonicity on the operator $F$ in a variational inequality problem (VIP), which is a particular case of MVIP (see Section II for more details), it is shown in [11], [12] that the solution of the VIP is globally asymptotically stable and globally exponentially stable, respectively, for the corresponding projected dynamical system. In contrast to the results mentioned above with asymptotic or exponential stability guarantees, a modified proximal dynamical system is introduced in this paper to guarantee convergence within a fixed time. In [13], the authors introduced the notion of finite-time stability of an equilibrium point, where the convergence of the solutions to the equilibrium point is guaranteed in a finite time. Under this notion, the settling time, or time of convergence, depends upon the initial conditions and can grow unbounded with the distance of the initial condition from an equilibrium point. A stronger notion, called fixed-time stability, is developed in [14], where the settling-time has a finite upper bound for all initial conditions.

While there is some work on finite- or fixed-time stable schemes for certain classes of convex optimization problems, to the best of the authors' knowledge, this is the first paper proposing fixed-time stable proximal dynamical systems for MVIPs or general non-smooth convex optimization problems. In [15], modified gradient flow schemes are introduced for unconstrained and constrained convex optimization problems, as well as for min-max problems posed as convex-concave optimization problems. The work in [15] only considers linear equality constraints and assumes that the objective function is continuously differentiable, satisfies strict convexity, or is gradient-dominated. In [16], the authors propose dynamical systems in the context of differential inclusions, such that any of the maximal Filippov solutions converge to a strict local minimizer of a given gradient-dominated objective function in a finite time (see also, [17]). The schemes proposed in this paper apply to a broader class of problems, namely, MVIPs, and smooth/non-smooth convex optimization problems arise as special cases of the general framework considered in this paper. The proposed work has two main contributions:

(i) A modified continuous-time proximal dynamical system for solving MVIPs is proposed, and the existence and uniqueness of solutions, as well as their convergence to the solution of the corresponding MVIP, are shown. Tools from fixed-time stability theory are leveraged to prove that trajectories of the modified proximal dynamical system converge to the solution of the associated MVIP within a fixed time, irrespective of the initial conditions. Thus, the results in [15] are a special case of the ones presented in this paper;
(ii) Inspired from the ideas presented in [18], it is shown that the solution obtained using the forward-Euler dis-

Kunal Garg is with the Department of Aerospace Engineering, University of Michigan, Ann Arbor, MI 48109, USA; kgarg@umich.edu. Mayank Baranwal is with the Data & Decision Sciences, Tata Consultancy Services Research, Mumbai, Maharashtra 400607, India; baranwal.mayank@tcs.com. Rohit Gupta was with the Department of Mechanical Engineering, University of Michigan, Ann Arbor, MI 48109, USA. He is currently with the Department of Aerospace Engineering, IIT Bombay, Powai, Mumbai, Maharashtra 400076, India; rohit@aero.iitb.ac.in. Mouhacine Benosman is with Mitsubishi Electric Research Laboratories, Cambridge, MA 02139, USA; benosman@merl.com.

cretization of a general class of differential inclusions, with a fixed-time stable equilibrium point, converges to an arbitrarily small neighborhood of it in a fixed number of time-steps, independent of the initial conditions. Prior works on fixed-time stability do not provide an analysis of the convergence behavior upon discretization. In contrast, our results bridge this gap and provide rigorous analysis of the convergence of the discretized dynamics.

In what follows, an inner product on $\mathbb{R}^n$ is denoted by $\langle\cdot,\cdot\rangle$, and $\|\cdot\| := \sqrt{\langle\cdot,\cdot\rangle}$ always denotes the induced norm except in some places, where it may denote an arbitrary norm on $\mathbb{R}^n$ when it is clear from the context.

## II. Problem Description

This paper considers the MVIP of the form:

$$\text{Find } x^* \in \mathbb{R}^n \text{ such that}$$
$$\langle F(x^*), x - x^*\rangle + g(x) - g(x^*) \geq 0 \text{ for all } x \in \mathbb{R}^n, \quad (1)$$

where $F : \operatorname{dom} g \to \mathbb{R}^n$ is an operator and $g : \mathbb{R}^n \to \mathbb{R} \cup \{\infty\}$ is a proper, lower semi-continuous convex function with $\operatorname{dom} g := \{x \in \mathbb{R}^n : g(x) < \infty\}$. Note that solving the MVIP is equivalent to the problem of solving a generalized equation of the form:

$$\text{Find } x^* \in \mathbb{R}^n \text{ such that } 0 \in F(x^*) + \partial g(x^*),$$

where the sub-differential mapping $\partial g : \mathbb{R}^n \rightrightarrows \mathbb{R}^n$ is a maximal monotone operator (see [19]). The function $g$ in (1) is not necessarily differentiable, for instance, $g$ could represent the indicator function of some non-empty, closed convex set $\mathcal{C} \subseteq \mathbb{R}^n$, i.e., $g = \delta_{\mathcal{C}}$, where

$$\delta_{\mathcal{C}}(x) = \begin{cases} 0, & \text{if } x \in \mathcal{C}; \\ \infty, & \text{otherwise,} \end{cases}$$

in which case, the MVIP reduces to a VIP. This paper is concerned with the following problem:

**Problem 1.** *Design a continuous-time proximal dynamical system, such that its solution converges to the solution of the MVIP (1) within a fixed time, independent of the initial conditions.*

**Remark 1.** *Note that the following optimization problem:*

$$\min_{x\in\mathbb{R}^n} f(x) + h(x),$$

*with $f : \operatorname{dom} h \to \mathbb{R}$ being a differentiable convex function and $h : \mathbb{R}^n \to \mathbb{R} \cup \{\infty\}$ being a proper, lower semi-continuous convex function, is equivalent to an MVIP with the operator $F = \nabla f$ and the function $g = h$ in (1).*

## III. Preliminaries

Consider the autonomous differential equation:

$$\dot{x} = X(x), \quad (2)$$

where the vector field $X : \mathbb{R}^n \to \mathbb{R}^n$ is continuous and $X(x^*) = 0$ for some $x^* \in \mathbb{R}^n$.[1]

[1] It is assumed that the solutions of (2) exist and are uniquely determined.

**Definition 1.** *The equilibrium point $x^*$ of (2) is said to be **fixed-time stable** if it is Lyapunov stable and $\lim_{t\to T(x(0))} x(t) = x^*$, where $\sup_{x(0)\in\mathbb{R}^n} T(x(0)) < \infty$ and $T : \mathbb{R}^n \to [0,\infty)$ is the settling-time function.*

**Lemma 1** ([14]). *Suppose that there exists a radially unbounded, continuously differentiable function $V : \mathbb{R}^n \to [0,\infty)$ such that $V(x^*) = 0, V(x) > 0$ for all $x \in \mathbb{R}^n \setminus \{x^*\}$ and the time-derivative of the function $V$ along the trajectories of (2) satisfies $\dot{V}(x) \leq -(a_1V(x)^{\gamma_1} + a_2V(x)^{\gamma_2})^{\gamma_3}$ for all $x \in \mathbb{R}^n \setminus \{x^*\}$, with $a_1, a_2, \gamma_1, \gamma_2, \gamma_3 > 0$ such that $\gamma_1\gamma_3 < 1$ and $\gamma_2\gamma_3 > 1$. Then, the equilibrium point $x^*$ of (2) is fixed-time stable with $T(x(0)) \leq \frac{1}{a_1{}^{\gamma_3}(1-\gamma_1\gamma_3)} + \frac{1}{a_2{}^{\gamma_3}(\gamma_2\gamma_3-1)}$, for any initial condition $x(0) \in \mathbb{R}^n$.*

**Remark 2.** *Lemma 1 provides characterization of fixed-time stability in terms of a Lyapunov function $V$. The existence of such a Lyapunov function for a suitably modified proximal dynamical system constitutes the foundation for the analysis in the paper, where Lemma 1 is used with $\gamma_3 = 1$.*

Next, some well-known definitions of the various notions of monotonicity of operators are given below (see, e.g., [8] for more details).

**Definition 2.** *An operator $F : \Omega \to \mathbb{R}^n$, where $\Omega$ is a non-empty subset of $\mathbb{R}^n$, is called:*

*(i)* ***Monotone****, if for all $x, y \in \Omega$,*
$$\langle F(x) - F(y), x - y\rangle \geq 0.$$
*(ii)* ***Strongly monotone*** *with modulus $\mu$, if there exists $\mu > 0$ such that for all $x, y \in \Omega$,*
$$\langle F(x) - F(y), x - y\rangle \geq \mu\|x - y\|^2.$$
*(iii)* ***Pseudomonotone****, if for all $x, y \in \Omega$,*
$$\langle F(y), x - y\rangle \geq 0 \implies \langle F(x), x - y\rangle \geq 0.$$
*(iv)* ***Strongly pseudomonotone*** *with modulus $\mu$, if there exists $\mu > 0$ such that for all $x, y \in \Omega$,*
$$\langle F(y), x - y\rangle \geq 0 \implies \langle F(x), x - y\rangle \geq \mu\|x - y\|^2.$$

The following lemma will be required in the proof of the first main result of the paper.

**Lemma 2.** *For any given $c \in (0,1)$, let $\varepsilon(c) = \frac{\log(c)}{\log\left(\frac{1-c}{1+c}\right)} > 0$. Then, the following strict inequality:*

$$\left(\frac{1-c}{1+c}\right)^{1-\alpha_1} > c,$$

*holds for each $\alpha_1 \in (1 - \varepsilon(c), 1]$. Furthermore, the following strict inequality:*

$$\left(\frac{1-c}{1+c}\right)^{\alpha_2-1} > c,$$

*holds for each $\alpha_2 \in [1, 1 + \varepsilon(c))$.*

The proof of the above lemma can be completed using simple algebraic manipulations and hence, it is omitted in the interest of space. The following assumptions will always be in place for the rest of this paper unless stated otherwise:

**Assumption 1.** *The operator $F$ is:*

*(i) Strongly monotone with modulus $\mu$.*
*(ii) Lipschitz continuous with Lipschitz constant $L$.*

## IV. Modified Proximal Dynamical System

Recall that the proximal operator associated with a proper, lower semi-continuous convex function $w : \mathbb{R}^n \to \mathbb{R} \cup \{\infty\}$ is defined as follows:

$$\operatorname{prox}_w(x) := \arg\min_{y\in\mathbb{R}^n} \left( w(y) + \frac{1}{2}\|x-y\|^2 \right). \tag{3}$$

For solving the MVIP (1), first consider the following nominal proximal dynamical system:

$$\dot{x} = -\kappa \left( x - \operatorname{prox}_{\lambda g}(x - \lambda F(x)) \right), \tag{4}$$

where $\kappa, \lambda > 0$ are some constants. In what follows, for the sake of brevity, set

$$y(x) := \operatorname{prox}_{\lambda g}(x - \lambda F(x)), \tag{5}$$

where $x \in \mathbb{R}^n$. The following lemma establishes the relationship between an equilibrium point of the nominal proximal dynamical system and a solution of the associated MVIP.[2]

**Lemma 3.** *A point $\bar{x} \in \mathbb{R}^n$ is an equilibrium point of (4) if and only if it solves (1).*

*Proof:* From [20, Proposition 12.26], it follows that

$$\begin{aligned} \bar{x} = y(\bar{x}) \iff& \langle (\bar{x} - \lambda F(\bar{x})) - \bar{x}, z - \bar{x}\rangle + \lambda g(\bar{x}) \le \lambda g(z), \\ \iff& \lambda \langle F(\bar{x}), z - \bar{x}\rangle + \lambda g(z) - \lambda g(\bar{x}) \ge 0, \\ \iff& \langle F(\bar{x}), z - \bar{x}\rangle + g(z) - g(\bar{x}) \ge 0, \end{aligned}$$

for all $z \in \mathbb{R}^n$. Hence, $\bar{x} \in \mathbb{R}^n$ is an equilibrium point of (4) if and only if it solves (1). ■

The following lemma will also be required to prove the first main result of the paper.

**Lemma 4.** *For each given $\lambda \in \left(0, \frac{2\mu}{L^2}\right)$, let $c = \sqrt{1 - 2\lambda\mu + \lambda^2 L^2} \in (0,1)$. Then, the inequality, $\|y(x) - x^*\| \le c\|x - x^*\|$, holds for all $x \in \mathbb{R}^n$, where the operator $y$ is defined as in (5) and $x^* \in \mathbb{R}^n$ is a solution of (1).*

*Proof:* The proof follows from the arguments used in the proof of [20, Proposition 26.16 (ii)], with $\alpha = \mu$, $\beta = L$, $\gamma = \lambda$ and $\tau = c$, however, for the sake of the readers' convenience, a complete proof is given here. Recall that for a proper, lower semi-continuous convex function $g$, the operator $y$ is firmly non-expansive (see, e.g., [20, Corollary 23.11 (i)]). Hence, the following inequality:

$$\begin{aligned} \|y(x) - y(x^*)\|^2 &\le \|(x - \lambda F(x)) - (x^* - \lambda F(x^*))\|^2 \\ &= \|x - x^*\|^2 - 2\lambda\langle F(x) - F(x^*), x - x^*\rangle \\ &\quad + \lambda^2\|F(x) - F(x^*)\|^2 \end{aligned} \tag{6}$$

holds for all $x \in \mathbb{R}^n$. From Lemma 3, it follows that $x^* \in \mathbb{R}^n$ is also an equilibrium point of (4) and hence, it follows that $y(x^*) = x^*$. Using this fact and Assumption 1, the right hand side of (6) can further be upper bounded and so, (6) results into:

$$\begin{aligned} \|y(x) - x^*\|^2 &\le \|x - x^*\|^2 - 2\lambda\mu\|x - x^*\|^2 \\ &\quad + \lambda^2 L^2\|x - x^*\|^2, \\ &= \underbrace{\left(1 - 2\lambda\mu + \lambda^2 L^2\right)}_{:=\bar{c}} \|x - x^*\|^2. \end{aligned} \tag{7}$$

Note that $\bar{c} \in (0,1)$, since by the assumption of the theorem, $\lambda \in \left(0, \frac{2\mu}{L^2}\right)$, and so, (7) can be re-written as $\|y(x) - x^*\| \le c\|x - x^*\|$, where $c := \sqrt{\bar{c}} \in (0,1)$, which completes the proof, since $c$ is independent of the choice of $x \in \mathbb{R}^n$. ■

A modified proximal dynamical system is now introduced, such that its equilibrium point is fixed-time stable. Let $\operatorname{Fix}(y) := \{\bar{x} \in \mathbb{R}^n : y(\bar{x}) = \bar{x}\}$ and consider the modification of (4) given by:

$$\dot{x} = -\rho(x)\left(x - y(x)\right), \tag{8}$$

where

$$\rho(x) := \begin{cases} 0, & \text{if } x \in \operatorname{Fix}(y); \\ \kappa_1 \frac{1}{\|x - y(x)\|^{1-\alpha_1}} + \kappa_2 \frac{1}{\|x - y(x)\|^{1-\alpha_2}}, & \text{otherwise}, \end{cases} \tag{9}$$

with $\kappa_1, \kappa_2 > 0$, $\alpha_1 \in (0,1)$ and $\alpha_2 > 1$. Note that in (9), the first term corresponding to the exponent $\alpha_1$ results in the finite-time stability of the equilibrium point of (8), while the second term corresponding to the exponent $\alpha_2$ helps in bounding the time of convergence to the equilibrium point of (8), uniformly for all initial conditions (also see [15, Remark 4]). The following lemma establishes the relationship between equilibrium points of the modified and nominal proximal dynamical systems.

**Lemma 5.** *A point $\bar{x} \in \mathbb{R}^n$ is an equilibrium point of (8) if and only if it is an equilibrium point of (4).*

*Proof:* Using (9), it is clear that if $\bar{x} \in \mathbb{R}^n$ is an equilibrium point of (8), i.e., $\bar{x} \in \operatorname{Fix}(y)$, then it is also an equilibrium point of (4). To show the other implication, it suffices to note that $\rho(x) = 0$ for any $x \in \operatorname{Fix}(y)$. ■

The following proposition establishes that the solutions of (8) exist and are uniquely determined for all forward times.

**Proposition 1.** *Let $X : \mathbb{R}^n \to \mathbb{R}^n$ be a locally Lipschitz continuous vector field such that*

$$X(\bar{x}) = 0 \text{ and } \langle x - \bar{x}, X(x)\rangle > 0$$

*for all $x \in \mathbb{R}^n \setminus \{\bar{x}\}$, where $\bar{x} \in \mathbb{R}^n$. Consider the following autonomous differential equation:*

$$\dot{x} = -\sigma(x)X(x), \tag{10}$$

*where*

$$\sigma(x) := \begin{cases} 0, & \text{if } X(x) = 0; \\ \kappa_1 \frac{1}{\|X(x)\|^{1-\alpha_1}} + \kappa_2 \frac{1}{\|X(x)\|^{1-\alpha_2}}, & \text{otherwise}, \end{cases}$$

*with $\kappa_1, \kappa_2 > 0$, $\alpha_1 \in (0,1)$ and $\alpha_2 > 1$. Then, the right-hand side of (10) is continuous for all $x \in \mathbb{R}^n$, and starting from any given initial condition, a solution of (10) exists and is uniquely determined for all $t \ge 0$.*

The proof of the above proposition is given in Appendix A.

**Remark 3.** *In the case, when the vector field $X$ is chosen to be the one in (4), i.e., $X(x) := x - y(x)$ for any $x \in \mathbb{R}^n$, then*

[2] In the special case, when $F = \nabla f$ in (1), where $f : \operatorname{dom} g \to \mathbb{R}$ is a differentiable convex function, Lemma 3 follows from [20, Corollary 27.3].

*it can be shown that the vector field $X$ has the property:*

$$\langle x-\bar{x}, X(x)\rangle > 0 \tag{11}$$

*for all $x \in \mathbb{R}^n \setminus \{\bar{x}\}$, where $\bar{x} \in \mathrm{Fix}(y)$. To see this, first note that from [21, Theorem 3.1] and Lemma 3, it follows that the vector field in (4) has a unique equilibrium point $\bar{x} = x^*$, where $x^* \in \mathbb{R}^n$ is the solution of (1), i.e., the set $\mathrm{Fix}(y)$ consists only of a single element $\bar{x} = x^*$. Furthermore, the following equality:*

$$\langle x-\bar{x}, x-y(x)\rangle = \|x-\bar{x}\|^2 + \langle x-\bar{x}, \bar{x}-y(x)\rangle, \tag{12}$$

*holds for all $x \in \mathbb{R}^n$. Using the Cauchy–Schwarz inequality and Lemma 4 (keeping in mind the fact that $\bar{x} = x^*$), the second term in the right hand side of (12) can be lower bounded and so, (12) results into:*

$$\langle x-\bar{x}, X(x)\rangle = \langle x-\bar{x}, x-y(x)\rangle \geq (1-c)\|x-\bar{x}\|^2,$$

*where $c \in (0,1)$, from which, it follows that (11) holds for all $x \in \mathbb{R}^n \setminus \{\bar{x}\}$.*

The following theorem establishes the first main result of the paper.

**Theorem 1.** *For any given $\lambda \in \left(0, \frac{2\mu}{L^2}\right)$, let $c = \sqrt{1-2\lambda\mu+\lambda^2L^2} \in (0,1)$ and $\varepsilon(c) = \frac{\log(c)}{\log\left(\frac{1-c}{1+c}\right)} > 0$. Then, the solution $x^* \in \mathbb{R}^n$ of (1) is a fixed-time stable equilibrium point of (8) for any $\alpha_1 \in (1-\varepsilon(c), 1) \cap (0,1)$ and $\alpha_2 \in (1, 1+\varepsilon(c))$.*

*Proof:* First note that the vector field in (4) is Lipschitz continuous on $\mathbb{R}^n$, which follows from the Lipschitz continuity of the proximal operator (see, e.g., [20, Proposition 12.28]) and Assumption 1, with a unique equilibrium point $\bar{x} = x^*$ (see Remark 3). Furthermore, it also satisfies the required properties assumed in Proposition 1 (see Remark 3). Hence, from Proposition 1, it follows that starting from any given initial condition, a solution of (8) exists and is uniquely determined for all forward times. Consider now a radially unbounded candidate Lyapunov function $V : \mathbb{R}^n \to [0,\infty)$ defined as $V(x) := \frac{1}{2}\|x-x^*\|^2$, where, from Lemma 5, $x^* \in \mathbb{R}^n$ is the unique equilibrium point of the vector field in (8). The time-derivative of $V$ along the solution of (8), starting from any $x(0) \in \mathbb{R}^n \setminus \{x^*\}$, reads:[3]

$$\begin{aligned}\dot{V} &= -\left\langle x-x^*, \kappa_1\frac{x-y(x)}{\|x-y(x)\|^{1-\alpha_1}} + \kappa_2\frac{x-y(x)}{\|x-y(x)\|^{1-\alpha_2}}\right\rangle \\ &= -\left\langle x-x^*, \kappa_1\frac{x-x^*}{\|x-y(x)\|^{1-\alpha_1}} + \kappa_2\frac{x-x^*}{\|x-y(x)\|^{1-\alpha_2}}\right\rangle \\ &\quad -\left\langle x-x^*, \kappa_1\frac{x^*-y(x)}{\|x-y(x)\|^{1-\alpha_1}} + \kappa_2\frac{x^*-y(x)}{\|x-y(x)\|^{1-\alpha_2}}\right\rangle.\end{aligned} \tag{13}$$

Using the Cauchy–Schwarz inequality, the second term in the right hand side of (13) can be bounded to obtain

$$\dot{V} \leq -\left(\kappa_1\frac{\|x-x^*\|^2}{\|x-y(x)\|^{1-\alpha_1}} + \kappa_2\frac{\|x-x^*\|^2}{\|x-y(x)\|^{1-\alpha_2}}\right)$$
$$+\left(\kappa_1\frac{\|x-x^*\|\|x^*-y(x)\|}{\|x-y(x)\|^{1-\alpha_1}} + \kappa_2\frac{\|x-x^*\|\|x^*-y(x)\|}{\|x-y(x)\|^{1-\alpha_2}}\right). \tag{14}$$

Note that by the assumption of the Theorem, $\lambda \in \left(0, \frac{2\mu}{L^2}\right)$ and so, Lemma 4 can be invoked. Using the triangle inequality and Lemma 4, there exists $c \in (0,1)$ such that the following inequality:

$$\|x-y(x)\| \leq \|x-x^*\|+\|y(x)-x^*\| \leq (1+c)\|x-x^*\| \tag{15}$$

holds for all $x \in \mathbb{R}^n$. Similarly, using the reverse triangle inequality and Lemma 4, there exists $c \in (0,1)$ such that the following inequality:

$$\|x-y(x)\| \geq \|x-x^*\|-\|y(x)-x^*\| \geq (1-c)\|x-x^*\| \tag{16}$$

also holds for all $x \in \mathbb{R}^n$. Using (15), (16) and Lemma 4, the right hand side of (14) can further be upper bounded and so, (14) results into:

$$\begin{aligned}\dot{V} \leq &-\left(\frac{\kappa_1}{(1+c)^{1-\alpha_1}}\frac{\|x-x^*\|^2}{\|x-x^*\|^{1-\alpha_1}} + \frac{\kappa_2}{(1-c)^{1-\alpha_2}}\frac{\|x-x^*\|^2}{\|x-x^*\|^{1-\alpha_2}}\right) \\ &+\left(\frac{c\kappa_1}{(1-c)^{1-\alpha_1}}\frac{\|x-x^*\|^2}{\|x-x^*\|^{1-\alpha_1}} + \frac{c\kappa_2}{(1+c)^{1-\alpha_2}}\frac{\|x-x^*\|^2}{\|x-x^*\|^{1-\alpha_2}}\right) \\ &= -q_1(\kappa_1,\alpha_1)\|x-x^*\|^{1+\alpha_1} - q_2(\kappa_2,\alpha_2)\|x-x^*\|^{1+\alpha_2},\end{aligned} \tag{17}$$

where $q_1(\kappa_1,\alpha_1) := \frac{\kappa_1}{(1-c)^{1-\alpha_1}}\left(\left(\frac{1-c}{1+c}\right)^{1-\alpha_1} - c\right)$ and $q_2(\kappa_2,\alpha_2) := \frac{\kappa_2}{(1+c)^{1-\alpha_2}}\left(\left(\frac{1-c}{1+c}\right)^{\alpha_2-1} - c\right)$. From Lemma 2, it follows that there exists $\varepsilon(c) = \frac{\log(c)}{\log\left(\frac{1-c}{1+c}\right)} > 0$ such that $q_1(\kappa_1,\alpha_1) > 0$ for any $\alpha_1 \in (1-\varepsilon(c),1) \cap (0,1)$, and $q_2(\kappa_2,\alpha_2) > 0$ for any $\alpha_2 \in (1, 1+\varepsilon(c))$. Hence, (17) results into:

$$\dot{V} \leq -\left(a_1(\kappa_1,\alpha_1)V^{\gamma(\alpha_1)} + a_2(\kappa_2,\alpha_2)V^{\gamma(\alpha_2)}\right), \tag{18}$$

where $a_i(\kappa_i,\alpha_i) := 2^{\gamma(\alpha_i)}q_i(\kappa_i,\alpha_i)$, $i=1,2$ and $\gamma(\alpha) := \frac{1+\alpha}{2}$. Note that $a_1(\kappa_1,\alpha_1) > 0$, $\gamma(\alpha_1) \in (0.5,1)$ for any $\alpha_1 \in (1-\varepsilon(c),1) \cap (0,1)$ and $a_2(\kappa_2,\alpha_2) > 0$, $\gamma(\alpha_2) > 1$ for any $\alpha_2 \in (1,1+\varepsilon(c))$. Hence, the proof can be concluded using Lemma 1. ■

In the special case, when the function $w$ in (3) is chosen to be the indicator function of a non-empty, closed convex set $\mathcal{C} \subseteq \mathbb{R}^n$, the proximal operator reduces to the projection operator, i.e., $\mathsf{P}_{\mathcal{C}} = \mathrm{prox}_{\delta_{\mathcal{C}}}$, where the projection operator is defined as $\mathsf{P}_{\mathcal{C}}(x) := \arg\min_{y\in\mathcal{C}} \frac{1}{2}\|x-y\|^2$, and so, the nominal proximal dynamical system reduces to a nominal projected dynamical system:

$$\dot{x} = -\kappa\left(x - \mathsf{P}_{\mathcal{C}}(x-\lambda F(x))\right), \tag{19}$$

with $\kappa, \lambda > 0$, which can be used to solve VIPs (see, e.g., [3], [8], [12]).[4] Furthermore, the modified proximal dynamical system now reduces to a modified projected dynamical system:

$$\dot{x} = -\rho(x)\left(x - \mathsf{P}_{\mathcal{C}}(x-\lambda F(x))\right). \tag{20}$$

It is shown in [3], [8], [12] that the equilibrium point of (19) is globally exponentially stable for a strongly pseudomonotone and Lipschitz continuous operator $F$. The following corollary of Theorem 1 establishes the fixed-time stability of the equilibrium point of the modified projected dynamical system (20).

[3]For the sake of brevity, the expressions $\dot{V}(x(t))$, $x(t)$ and $y(x(t))$ are abbreviated as $\dot{V}$, $x$ and $y(x)$, respectively, in the proof.

[4]The existence and uniqueness of the solution of VIP holds for a strongly pseudomonotone and Lipschitz continuous operator $F$ ([22, Theorem 2.1]).

**Corollary 1.** *For any given $\lambda \in \left(0, \frac{2\mu}{L^2}\right)$, there exist $c \in (0,1)$ and $\varepsilon(c) > 0$ as given in Theorem 1, such that the solution $x^* \in \mathcal{C}$ of (1), with $g = \delta_{\mathcal{C}}$, where $\mathcal{C}$ is a closed convex set, is a fixed-time stable equilibrium point of (20) for any $\alpha_1 \in (1-\varepsilon(c), 1) \cap (0,1)$ and $\alpha_2 \in (1, 1+\varepsilon(c))$.*

**Remark 4.** *In the special case of a projection operator, Lemma 4, and hence, Corollary 1, continues to hold, even when (i) in Assumption 1 is relaxed to strong pseudomonotonicity (see the proof of [8, Theorem 2]). Furthermore, by following the steps given in the proof of Theorem 1, with $\kappa_2 = 0$ and $\alpha_1 = 1$, it can be seen that [8, Theorem 2] is now a special case of Corollary 1, from which only the exponential stability (instead of fixed-time stability) of the equilibrium point can be concluded.*

## V. $(T,\epsilon)$-Close Discrete-Time Approximation Scheme

Continuous-time dynamical systems, such as the one given by (8) offer effective insights into designing accelerated schemes for solving MVIPs. However, in practice, a discrete-time implementation is used for solving such problems using iterative methods. In general, the fixed-time convergence behavior of the continuous-time dynamical system might not be preserved in the discrete-time setting. A consistent discrete-time approximation scheme preserves the convergence behavior of the continuous-time dynamical system in the discrete-time setting. The authors in [23] present a consistent *semi-implicit* discrete-time approximation scheme (called as *consistent discrete-time approximation*) for practically fixed-time stable systems. In contrast, inspired from the ideas presented in [18], this section characterizes sufficient conditions that lead to an *explicit* $(T,\epsilon)$-close discrete-time approximation scheme for the fixed-time convergent modified proximal dynamical system. In order to achieve this goal, the following adaptation of [24, Definition 3.2] is needed:

**Definition 3.** *Let $T, \epsilon, \eta > 0$ be given and consider a solution $x_c : [0,T] \to \mathbb{R}^n$ of the following autonomous differential equation:*

$$\dot{x} = X_c(x), \;\; x(0) = x_{c,0},$$

*where $X_c : \mathbb{R}^n \to \mathbb{R}^n$ is a continuous vector field, and a solution $x_d : \left\{0, 1, \ldots, \left\lfloor \frac{T}{\eta} \right\rfloor\right\} \to \mathbb{R}^n$ of the following autonomous difference equation:*

$$x_{k+1} = X_d(x_k), \;\; x_0 = x_{d,0},$$

*where $X_d : \mathbb{R}^n \to \mathbb{R}^n$ is a continuous map. The solutions $x_c$ and $x_d$ are said to be $(T,\epsilon)$-close, if:*

*(i) For each $t \in [0,T]$, there exists $k \in \left\{0, 1, \ldots, \left\lfloor \frac{T}{\eta} \right\rfloor\right\}$ such that $|t - \eta k| < \epsilon$ and $\|x_c(t) - x_d(k)\| < \epsilon$;*

*(ii) For each $k \in \left\{0, 1, \ldots, \left\lfloor \frac{T}{\eta} \right\rfloor\right\}$, there exists $t \in [0,T]$ such that $|t - \eta k| < \epsilon$ and $\|x_d(k) - x_c(t)\| < \epsilon$.*

Now, consider the forward-Euler discretization of (2):

$$x_{k+1} = x_k + \eta X(x_k), \tag{21}$$

where $\eta > 0$ is the time-step. The following theorem characterizes sufficient conditions that lead to an $(T,\epsilon)$-close discrete-time approximation scheme for a differential equation with a fixed-time stable equilibrium point.

**Theorem 2.** *Assume that the conditions of Lemma 1 hold, with $\gamma_1 = 1 - \frac{1}{\xi}, \gamma_2 = 1 + \frac{1}{\xi}$ and $\gamma_3 = 1$, where $\xi > 1$. Furthermore, assume that the function $V$ also satisfies the following lower bound for each $x \in \mathbb{R}^n$:*

$$V(x) \geq \beta(\|x - \bar{x}\|), \tag{22}$$

*where $\beta$ is a class-$\mathcal{K}_\infty$ function and $\bar{x}$ is the equilibrium point of (2).[5] Then, for any given $x_0 \in \mathbb{R}^n$ and $\epsilon > 0$, there exists $\eta^* > 0$ such that for any $\eta \in (0, \eta^*]$, the following holds:*

$$\|x_k - \bar{x}\| < \begin{cases} \beta^{-1}\left(\left(\sqrt{\frac{a_1}{a_2}} \tan\left(\frac{\pi}{2} - \frac{\sqrt{a_1 a_2}}{\xi}\eta k\right)\right)^{\xi}\right) + \epsilon, & k \leq \left\lceil \frac{\xi\pi}{2\eta\sqrt{a_1 a_2}} \right\rceil; \\ \epsilon, & \text{otherwise}, \end{cases} \tag{23}$$

*where $x_k$ is a solution of (21) starting from the point $x_0$.*

*Proof:* Under the conditions of Lemma 1, it follows that the following inequality:

$$V(x(t)) \leq \left(\sqrt{\frac{a_1}{a_2}} \tan\left(\tan^{-1}\left(\sqrt{\frac{a_2}{a_1}} V(x(0))^{\frac{1}{\xi}}\right) - \frac{\sqrt{a_1 a_2}}{\xi} t\right)\right)^{\xi} \tag{24}$$

holds for each $t \in [0, \hat{\tau}]$, with

$$\hat{\tau} = \frac{\xi}{\sqrt{a_1 a_2}} \tan^{-1}\left(\sqrt{\frac{a_2}{a_1}} V(x(0))^{\frac{1}{\xi}}\right).$$

Hence, using (22), it follows from (24) that the following holds:

$$\|x(t) - \bar{x}\| \leq \begin{cases} \beta^{-1}\left(\left(\sqrt{\frac{a_1}{a_2}} \tan\left(\frac{\pi}{2} - \frac{\sqrt{a_1 a_2}}{\xi} t\right)\right)^{\xi}\right), & t < \frac{\xi\pi}{2\sqrt{a_1 a_2}}; \\ 0, & \text{otherwise}, \end{cases} \tag{25}$$

where $x : [0,\infty) \to \mathbb{R}^n$ is a solution of (2) starting from the initial condition $x_0 \in \mathbb{R}^n$. Furthermore, it can be verified that all the requirements as stated in [24, Theorem 5.2] are met. Hence, for each $\epsilon > 0$ and each $T \geq 0$, there exists $\eta^* > 0$ with the following property: for any $\eta \in (0, \eta^*]$ and a solution $x_k$ of (21) starting from the point $x_0$, there exists a solution $x$ of (2) starting from the point $x_0$ such that the solutions $x$ and $x_k$ are $(T,\epsilon)$-close.

Note that from the triangle inequality, it follows that for any given $k \in \{0, 1, \ldots\}$, the following inequality:

$$\|x_k - \bar{x}\| \leq \|x(t) - \bar{x}\| + \|x_k - x(t)\| \tag{26}$$

holds for each $t \in [0, \infty)$. For any given $\eta \in (0, \eta^*]$, substituting now $t = \eta k$ in (26) and then using (25) and the $(T,\epsilon)$-closeness of the solutions $x_k$ and $x$, yields (23), which completes the proof. ∎

[5] A function $\beta : [0,\infty) \to [0,\infty)$ is of class-$\mathcal{K}_\infty$, if it is continuous and strictly increasing, with $\beta(0) = 0$ and $\lim_{r\to\infty} \beta(r) = \infty$.

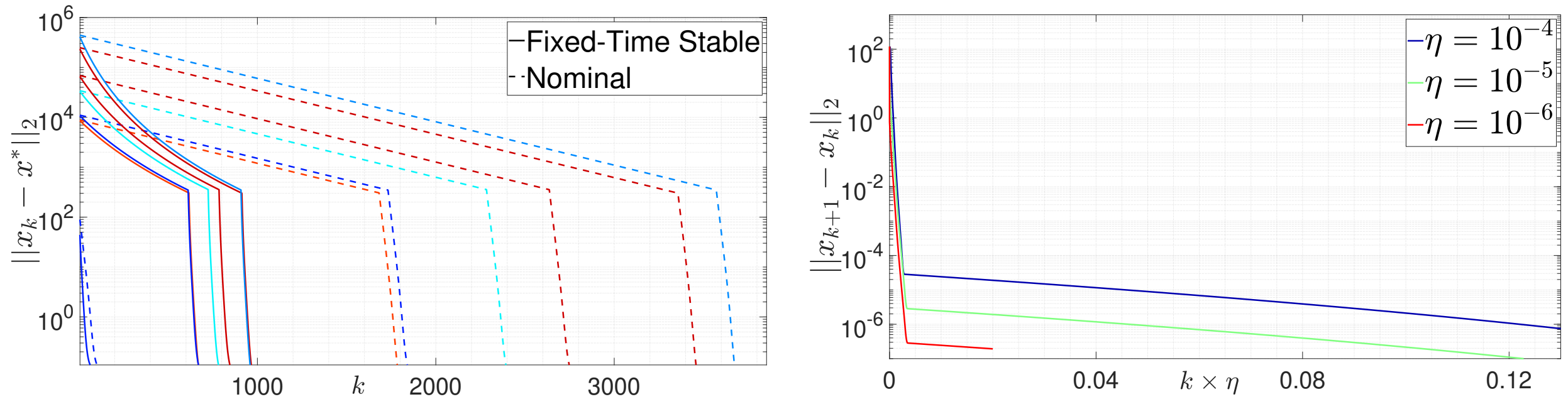

(a) (b)

Fig. 1. Plot of $\|x_k - x^*\|_2$ vs. the number of iterations $k$ for various initial conditions $x(0) \in \mathbb{R}^3$ in (a) and of $\|x_{k+1} - x_k\|_2$ for various values of the time-step $\eta$ in (b).

The following result states that a discrete-time system obtained through the forward-Euler discretization of (8), is an $(T, \epsilon)$-close discrete-time approximation scheme, and the resulting discrete-time trajectories reach an $\epsilon$-neighborhood of the solution of the associated MVIP within a fixed number of steps, independent of the initial conditions.

**Corollary 2.** *Consider the forward-Euler discretization of the modified proximal dynamical system* (8)*:*

$$x_{k+1} = x_k - \eta\rho(x_k)\left(x_k - y(x_k)\right), \tag{27}$$

*where* $\eta > 0$ *is the time-step and* $\rho$ *is given by* (9)*, with* $\kappa_1, \kappa_2 > 0$*,* $\alpha_1(\xi) = 1 - \frac{2}{\xi}$ *and* $\alpha_2(\xi) = 1 + \frac{2}{\xi}$*, where* $\xi > \max\left\{2, \frac{2}{\varepsilon(c)}\right\}$*, with* $\varepsilon(c)$ *being as given in Theorem 1. Furthermore, for any given* $x_0 \in \mathbb{R}^n$*,* $\epsilon > 0$ *and* $\lambda \in \left(0, \frac{2\mu}{L^2}\right)$*, let* $a_1 = q_1(\kappa_1, \alpha_1(\xi))$*,* $a_2 = q_2(\kappa_2, \alpha_2(\xi))$*, where* $c$ *is as given in Lemma 4,* $\varepsilon$ *is as given in Lemma 2 and the functions* $q_1$*,* $q_2$ *are as defined in the proof of Theorem 1. Then, there exists* $\eta^* > 0$ *such that for any* $\eta \in (0, \eta^*]$*, the following holds:*

$$\|x_k - x^*\| < \begin{cases} \sqrt{2}\left(\sqrt{\frac{a_1}{a_2}}\tan\left(\frac{\pi}{2} - \frac{\sqrt{a_1 a_2}}{\xi}\eta k\right)\right)^{\frac{\xi}{2}} + \epsilon, & k \le k^*; \\ \epsilon, & \textit{otherwise}, \end{cases}$$

*where* $k^* = \left\lceil \frac{\xi\pi}{2\eta\sqrt{a_1 a_2}} \right\rceil$ *and* $x_k$ *is a solution of* (27) *starting from the point* $x_0$ *and* $x^* \in \mathbb{R}^n$ *is the solution of* (1)*.*

*Proof:* First observe that from the proof of Theorem 1, it follows that for any given $\lambda \in \left(0, \frac{2\mu}{L^2}\right)$, (18) holds as long as $\alpha_1(\xi) \in (1 - \varepsilon(c), 1) \cap (0, 1)$, with $\varepsilon(c) = \frac{\log(c)}{\log\left(\frac{1-c}{1+c}\right)} > 0$ and $\alpha_2(\xi) > 1$. Note that this former requirement means that $\xi > \max\left\{2, \frac{2}{\varepsilon(c)}\right\}$, while the latter one is always satisfied for any choice of $\xi > 2$. Hence, for any given $\lambda \in \left(0, \frac{2\mu}{L^2}\right)$, (18) holds for any $\xi > \max\left\{2, \frac{2}{\varepsilon(c)}\right\}$, with $q_1(\kappa_1, \alpha_1(\xi)) > 0$, $q_2(\kappa_2, \alpha_2(\xi)) > 0$, $\gamma(\alpha_1(\xi)) = 1 - \frac{1}{\xi}$ and $\gamma(\alpha_2(\xi)) = 1 + \frac{1}{\xi}$. Hence, all the requirements stated in Theorem 2 are met with the Lyapunov function $V(x) = \frac{1}{2}\|x - x^*\|^2$, $a = a_1$ and $b = a_2$. The proof now follows by invoking Theorem 2. ■

It is illustrated in [15], [16], [25] through numerical experiments that the solution obtained using the forward-Euler discretization of finite- and fixed-time stable gradient flows exhibits a superior rate of convergence as compared to the nominal gradient descent method. In contrast, the above corollary shows this analytically that for any given $\epsilon > 0$, the solution obtained using the forward-Euler discretization of (8), reaches an $\epsilon$-neighborhood of the solution of the associated MVIP within a fixed number of time steps, independent of the initial conditions.

## VI. Numerical Example

The fixed-time convergent behavior of the modified proximal dynamical system is illustrated through an example. The simulations are performed in the `R2022a` version of `MATLAB` on a desktop machine with `16GB DDR4 RAM` and an `Intel Core i7-7500` processor with a speed of `2.7 GHz`. The results are shown in `log-lin` plots for better visualization. Consider the following elastic-net logistic regression problem:

$$\min_{x \in \mathbb{R}^3} \sum_{i=1}^{100} \log\left(1 + \exp(-a_i b_i^\mathsf{T} x)\right) + \beta_1\|x\|_1 + \beta_2\|x\|_2^2,$$

where $a_i \in \{-1, 1\}$, $b_i \in \mathbb{R}^3$, $i = 1, \ldots, 100$ are chosen randomly using the `rand` command in `MATLAB`, $\beta_1, \beta_2 > 0$ and the non-smooth $\ell_1$-regularization term is added to prevent overfitting on the given data. The above elastic-net logistic regression problem can be solved in the framework of the MVIP (1), where the non-smooth $\ell_1$-penalty term plays the role of $g$ in (1). Interestingly, the proximal operator associated with it has a closed-form expression given by:

$$\text{prox}_{\upsilon\|\cdot\|_1}(x) = \text{sgn}(x)\max\{0, |x| - \upsilon\},$$

where $\upsilon > 0$. The following parameter set is chosen for numerical experimentation purposes: $\beta_1 = 2.5$, $\beta_2 = 0.25$, $\eta = 10^{-4}$, $\kappa_1 = 10$ and $\kappa_2 = 10$. It can be verified that Assumption 1 holds for this example (with $\mu = 2\beta_2 = 0.5$ and $L = 0.5337$) and hence, from Theorem 1, it follows that for any $\lambda \in (0, 3.51)$, the solution of the MVIP (1) is a fixed-time stable equilibrium point of (8).[6] Furthermore, the parameter $\lambda$ is chosen as $\lambda = 2.5$, which results in $c = 0.5294$ and $\varepsilon(c) = 0.5397$, and consequently, the parameters $\alpha_1, \alpha_2$ are chosen as $\alpha_1 = 0.5142$ and $\alpha_2 = 1.4858$.

[6]Note that $L = 2\beta + \frac{1}{4}\Lambda_{\max}\left([b_1 \cdots b_{100}][b_1 \cdots b_{100}]^\mathsf{T}\right)$, where $\Lambda_{\max}(A)$ denotes the maximum eigenvalue for a symmetric matrix $A$.

| Method | Run-Time (sec) | Function Evaluations (#) | Iterations (#) |
|---|---|---|---|
| `fmincon` | 0.014 | 200 | 30 |
| (27) | 0.01 | 19 | 19 |

TABLE I
COMPARISON BETWEEN FMINCON AND THE PROPOSED METHOD.

Figure 1(a) shows the distances between the solutions obtained using the forward-Euler discretization of (8) and the solution $x^*$ of (1), for various randomly chosen initial conditions, where $x^* \in \mathbb{R}^3$ is obtained using the `fmincon` function in `MATLAB`.[7] Solid lines show the convergence behavior of the forward-Euler discretization of the fixed-time stable proximal dynamical system (8) (labeled as 'Fixed-Time Stable'), while the dotted lines show the convergence behavior of the forward-Euler discretization of the nominal proximal dynamical system (4) (labeled as 'Nominal'). It can be observed that the forward-Euler discretization of the fixed-time stable proximal dynamical system (8) takes a significantly lesser number of iterations to converge to the same final error as compared to the forward-Euler discretization of the nominal proximal dynamical system (4), for various initial conditions $x(0) \in \mathbb{R}^3$ and it converges to a neighborhood of the solution in a fixed number of iterations. The theoretical upper bound on the settling-time function is found to be approximately 5.06 units, using which it is found that $k^* = 5.06 \times 10^4$.

The effect of varying the time-step $\eta \in \{10^{-4}, 10^{-5}, 10^{-6}\}$, is also studied for this example and the results are shown in Figure 1(b), where the convergence behavior is depicted as a function of the product of the number of iterates $k$ and $\eta$. As can be seen in Figure 1(b), below a suitably chosen value of the time-step $\eta$, the convergence of the iterates to a small error tolerance of approximately $10^{-4}$, is largely independent of $k \times \eta$. Hence, for smaller values of the time-step $\eta$, the number of iterations needed for convergence to a specified error tolerance is larger, however, having a smaller value of $\eta$ facilitates smaller steady-state error tolerance, since the discrete-time trajectory better approximates the continuous-time one. In summary, suitably chosen large values of the time-step $\eta$ enable faster convergence, i.e., it takes a fewer number of iterations to reach an $\epsilon$-neighborhood of the solution, albeit at the expense of a relatively larger steady-state tolerance $\epsilon$.

The convergence behavior of the forward-Euler discretization of the fixed-time stable proximal dynamical system (8) (given by (27)) is also compared with the `fmincon` function (using the active set algorithm) in `MATLAB`. It must be noted that the *built-in* `fmincon` function is heavily optimized for implementation (e.g., back-end routines implemented in `C`, parallel programming support, etc.), while the proposed method (27) has been implemented using a simple *loop-based* implementation in `MATLAB`. Despite the differences in the implementation of the two methods, Table I suggests that the `fmincon` function requires approximately 10 times more function evaluations as compared to the proposed method (27) for 1000 runs of this example, with random initializations and for the same optimality tolerance of $10^{-3}$. Recall that unlike the method proposed in this work, the `fmincon` function requires several function evaluations during each iteration as well. The computational times for the two methods are nearly the same, which might largely be due to the optimized implementation of the `fmincon` function in `MATLAB`.

[7]In Figure 1, $\|\cdot\|_2$ denotes the Euclidean norm.

## VII. CONCLUSIONS

In this paper, a modified proximal dynamical system is presented such that its solution exists, is uniquely determined, and converges to the unique solution of the associated MVIP in a fixed time, under the assumptions of strong monotonicity and Lipschitz continuity on the associated operator. Furthermore, as a special case for solving variational inequality problems, the proposed modified proximal dynamical system reduces to a fixed-time stable projected dynamical system, where the fixed-time stability of the modified projected dynamical system continues to hold, even if the assumption of strong monotonicity is relaxed to that of strong pseudomonotonicity. Finally, it is shown that the forward-Euler discretization of the modified proximal dynamical system results in an explicit $(T, \epsilon)$-close discrete-time approximation scheme.

## REFERENCES

[1] F. Facchinei and J.-S. Pang, *Finite-Dimensional Variational Inequalities and Complementarity Problems*. Springer, 2003.

[2] F. Giannessi, A. Maugeri, and P. M. Pardalos, *Equilibrium Problems: Nonsmooth Optimization and Variational Inequality Models*. Springer, 2001.

[3] E. Cavazzuti, M. Pappalardo, and M. Passacantando, "Nash equilibria, variational inequalities, and dynamical systems," *Journal of Optimization Theory and Applications*, vol. 114, no. 3, pp. 491–506, 2002.

[4] G. Scutari, D. P. Palomar, F. Facchinei, and J.-S. Pang, "Convex optimization, game theory, and variational inequality theory," *IEEE Signal Processing Magazine*, vol. 27, no. 3, pp. 35–49, 2010.

[5] V. Barbu, P. Neittaanmäki, and A. Niemistö, "Approximating optimal control problems governed by variational inequalities," *Numerical Functional Analysis and Optimization*, vol. 15, no. 5-6, pp. 489–502, 1994.

[6] P. Neittaanmäki and D. Tiba, "A variational inequality approach to constrained control problems for parabolic equations," *Applied Mathematics and Optimization*, vol. 17, no. 1, pp. 185–201, 1988.

[7] D. Kinderlehrer and G. Stampacchia, *An Introduction to Variational Inequalities and Their Applications*. SIAM, 2000.

[8] N. T. T. Ha, J.-J. Strodiot, and P. T. Vuong, "On the global exponential stability of a projected dynamical system for strongly pseudomonotone variational inequalities," *Optimization Letters*, vol. 12, no. 7, pp. 1625–1638, 2018.

[9] S. Hassan-Moghaddam and M. R. Jovanović, "Proximal gradient flow and Douglas–Rachford splitting dynamics: Global exponential stability via integral quadratic constraints," *Automatica*, vol. 123, 2021, article 109311.

[10] X. Hu and J. Wang, "Solving pseudomonotone variational inequalities and pseudoconvex optimization problems using the projection neural network," *IEEE Transactions on Neural Networks*, vol. 17, no. 6, pp. 1487–1499, 2006.

[11] Y. Xia and J. Wang, "A general methodology for designing globally convergent optimization neural networks," *IEEE Transactions on Neural Networks*, vol. 9, no. 6, pp. 1331–1343, 1998.

[12] Y. Xia, H. Leung, and J. Wang, "A projection neural network and its application to constrained optimization problems," *IEEE Transactions on Circuits and Systems I: Fundamental Theory and Applications*, vol. 49, no. 4, pp. 447–458, 2002.

[13] S. P. Bhat and D. S. Bernstein, "Finite-Time stability of continuous autonomous systems," *SIAM Journal on Control and Optimization*, vol. 38, no. 3, pp. 751–766, 2000.

[14] A. Polyakov, "Nonlinear feedback design for fixed-time stabilization of linear control systems," *IEEE Transactions on Automatic Control*, vol. 57, no. 8, pp. 2106–2110, 2012.

[15] K. Garg and D. Panagou, "Fixed-Time stable gradient flows: Applications to continuous-time optimization," *IEEE Transactions on Automatic Control*, vol. 66, no. 5, pp. 2002–2015, 2021.

[16] O. Romero and M. Benosman, "Finite-Time convergence in continuous-time optimization," in *Proceedings of the International Conference on Machine Learning*, 2020, pp. 8200–8209.
[17] ——, "Time-Varying continuous-time optimisation with pre-defined finite-time stability," *International Journal of Control*, pp. 1–18, 2020.
[18] S. Zhang, M. Benosman, O. Romero, and A. Cherian, "First-Order optimization inspired from finite-time convergent flows," *arXiv preprint arXiv:2010.02990v3*, 2021.
[19] R. T. Rockafellar, "On the maximal monotonicity of subdifferential mappings," *Pacific Journal of Mathematics*, vol. 33, no. 1, pp. 209–216, 1970.
[20] H. H. Bauschke and P. L. Combettes, *Convex Analysis and Monotone Operator Theory in Hilbert Spaces*. Springer, 2017.
[21] M. Aslam Noor, "Mixed variational inequalities," *Applied Mathematics Letters*, vol. 3, no. 2, pp. 73–75, 1990.
[22] P. T. Vuong and P. D. Khanh, "Qualitative properties of strongly pseudomonotone variational inequalities," *Optimization Letters*, vol. 10, no. 8, pp. 1669–1679, 2016.
[23] A. Polyakov, D. Efimov, and B. Brogliato, "Consistent discretization of finite-time and fixed-time stable systems," *SIAM Journal on Control and Optimization*, vol. 57, no. 1, pp. 78–103, 2019.
[24] R. G. Sanfelice and A. R. Teel, "Dynamical properties of hybrid systems simulators," *Automatica*, vol. 46, no. 2, pp. 239–248, 2010.
[25] P. Budhraja, M. Baranwal, K. Garg, and A. Hota, "Breaking the Convergence Barrier: Optimization via Fixed-Time Convergent Flows," in *Proceedings of the AAAI Conference on Artificial Intelligence*, vol. 36, no. 6, 2022, pp. 6115–6122.
[26] J. K. Hale, *Ordinary Differential Equations*. Krieger Publishing Company, 1980.
[27] C. Chicone, *Ordinary Differential Equations with Applications*. Springer, 2006.
[28] R. P. Agarwal and V. Lakshmikantham, *Uniqueness and Nonuniqueness Criteria for Ordinary Differential Equations*. World Scientific Publishing Company, 1993.

## Appendix A
## Proof of Proposition 1

*Proof:* First it is shown that the vector field in (10) is continuous on $\mathbb{R}^n$. To see this, it would suffice to note that the vector field in (10) is continuous at $\bar{x} \in \mathbb{R}^n$, since the function $\sigma$ is continuous on $\mathbb{R}^n \setminus \{\bar{x}\}$ and the vector field $X$ is assumed to be locally Lipschitz continuous on $\mathbb{R}^n$. To this end, note that $\lim_{x\to\bar{x}} \sigma(x)X(x) = 0$, since $\alpha_1 \in (0,1)$ and $\alpha_2 > 1$. The proof for the existence of a solution of (10), starting from any given initial condition for all forward times, is shown as follows. First note that the equilibrium point $\bar{x} \in \mathbb{R}^n$ of the vector field $X$ is unique since under its properties assumed in the proposition, it can be shown that the equilibrium point $\bar{x} \in \mathbb{R}^n$ of the vector field $-X$ is globally asymptotically stable and hence, it is unique (another way to see this, is to use the Cauchy–Schwarz inequality to upper bound the left-hand side of (11)). Using the fact that the vector field in (10) is continuous, it follows from [26, Theorem I.1.1] that for any given $x(0) \in \mathbb{R}^n$, there exists a solution of (10) on some interval $[0, \tau(x(0))]$, with $\tau(x(0)) > 0$. Furthermore, from [26, Theorem I.2.1] any such solution of (10) on the interval $[0, \tau(x(0))]$ has a continuation to a maximal interval of existence $[0, \bar{\tau}(x(0)))$. Consider now the candidate Lyapunov function $V : \mathbb{R}^n \to [0, \infty)$ defined as $V(x) := \frac{1}{2}\|x-\bar{x}\|^2$. The time-derivative of the candidate Lyapunov function $V$ along a solution of (10), starting from $x(0) \in \mathbb{R}^n$, reads:

$$\dot{V}(x(t)) = -\langle x(t) - \bar{x}, \sigma(x(t))X(x(t))\rangle . \tag{28}$$

Recalling that $\langle x - \bar{x}, X(x)\rangle \geq 0$ and $\sigma(x) \geq 0$ for all $x \in \mathbb{R}^n$, (28) results into:

$$\dot{V}(x(t)) \leq 0,$$

for all $t \in [0, \bar{\tau}(x(0)))$. Hence, $V(x(t)) \leq V(x(0))$ for all $t \in [0, \bar{\tau}(x(0)))$ and it follows that a solution of (10) defined on the interval $[0, \bar{\tau}(x(0)))$ lies entirely in the set $K_{x(0)} := \{z \in \mathbb{R}^n : \|z - \bar{x}\| \leq \|x(0) - \bar{x}\|\}$. Since the set $K_{x(0)}$ is compact, from [13, Proposition 2.1], it follows that $\bar{\tau}(x(0)) = \infty$. This completes the proof for the claim of the existence of the solution.

The proof for the uniqueness of the solution of (10), starting from any given initial condition for all forward times, is shown next. For any given $x(0) \in \mathbb{R}^n$, let $\gamma$ be a solution of (10), with $\gamma(0) = x(0)$ and consider the following two cases:

(i) In the first case, let $\gamma(0) \in \mathbb{R}^n \setminus \{\bar{x}\}$. It will be shown that a solution corresponding to the vector field in (10) is also a solution corresponding to the vector field $-X$, under a suitable reparameterization of time (see, e.g., [27, Section 1.5]). Let $\mathsf{T} := \inf\{t \geq 0 : \gamma(t) = \bar{x}\}$ and from the continuity of $\gamma$, it follows that $\mathsf{T} > 0$. Consider now the function $\mathsf{s} : [0, \mathsf{T}) \to [0, \infty)$ defined as follows:

$$\mathsf{s}(t) := \int_0^t \sigma(\gamma(\nu))d\nu. \tag{29}$$

Since the function $\sigma$ is continuous on $\mathbb{R}^n$, $\gamma$ is continuous on the interval $[0, \mathsf{T})$ and $\sigma(\gamma(\nu)) > 0$ for any $\nu \in [0, \mathsf{T})$, it follows that the function $\mathsf{s}$ is a strictly increasing continuous function, with $\frac{d\mathsf{s}}{dt} \neq 0$ for all $t \in (0, \mathsf{T})$. Furthermore, from the inverse function Theorem, it follows that the function $\mathsf{t} := \mathsf{s}^{-1}$ exists, is strictly increasing, continuous, and satisfies:

$$\left.\frac{d\mathsf{t}}{ds}\right|_{s=\mathsf{s}(t)} = \frac{1}{\sigma(\gamma(t))} \tag{30}$$

for all $t \in (0, \mathsf{T})$. Let $\bar{\gamma}(s) := \gamma(\mathsf{t}(s))$ and from the chain rule, it follows that

$$\frac{d\bar{\gamma}}{ds} = \left.\frac{d\gamma}{dt}\right|_{t=\mathsf{t}(s)} \frac{d\mathsf{t}}{ds}. \tag{31}$$

Using (30), (31) reads $\frac{d\bar{\gamma}}{ds} = -X(\bar{\gamma}(s))$. Hence, a solution corresponding to the vector field in (10) is also a solution corresponding to the vector field $X$, under the reparameterization of time given in (29). Furthermore, by following the steps, similar to the ones given in the proof of the first claim and recalling that the vector field $X$ is locally Lipschitz continuous on $\mathbb{R}^n$, it can be shown that for any given initial condition, there exists a unique solution corresponding to the vector field $X$ for all forward times (see, e.g., [28, Theorem 3.2.2]). Hence, $\bar{\gamma}$ is uniquely determined and since the function $\mathsf{s}$ is injective, with $\mathsf{s}(0) = 0$, it follows that $\gamma$ is also uniquely determined.

(ii) In the second case, let $\gamma(0) = \bar{x}$. Consider now the same candidate Lyapunov function $V$ as the one given in the first claim. By following the steps given in the proof of the first claim, it can be shown that the time-derivative of the candidate Lyapunov function $V$ along a solution of (10), starting from any given initial condition is always non-positive. From [28, Theorem 3.15.1], it follows that $\gamma$ is uniquely determined.

This completes the proof for the second claim. ■